\documentclass[12pt]{article}
\usepackage[parfill]{parskip}    
\usepackage{amsthm,amsmath,comment}
\usepackage{amssymb}
\usepackage[numbers]{natbib}
\usepackage{graphicx}
\def\vp{\psi}
\def\t{\theta}
\def\ve{\varepsilon}
\def\t{u}
\def\pa{\partial}
\def\a{\alpha}
\def\b{\beta}
\def\d{\delta}
\def\g{\gamma}
\def\G{\Gamma}
\def\C{{\cal C}}
\def\wm{\tilde{\mu}}
\def\p{p^\ast}
\def\Z{{\cal Z}}
\def\noi{\noindent}
\def\non{\nonumber}

\newtheorem{lem}{Lemma}
\newtheorem{thm}{Theorem}
\theoremstyle{definition}
\newtheorem{defn}{Definition}[section]
\title{Limit Theorems   for  Hybridization  Reactions on Oligonucleotide Microarrays}
\author{Grzegorz A. Rempala\thanks{Corresponding author. Department of Mathematics and Center for Genetics and Molecular Medicine, University of Louisville. } \  and Iwona Pawlikowska\thanks{Department of Mathematics, University of Silesia,  Katowice, Poland.}}
\date{Jan 1, 2008 \\ \ \\   }                                        
\begin{document}

\maketitle
\begin{abstract}
We derive herein the  limiting laws for certain  stationary distributions of birth-and-death processes related to  the classical model of  chemical adsorption-desorption reactions due to  Langmuir.   The model has been recently considered in the context of a  hybridization reaction  on an oligonucleotide DNA microarray. Our results imply   that   the truncated gamma- and beta- type  distributions can be  used as approximations to the  observed distributions of the  fluorescence readings of the oligo-probes on a microarray.  These findings might be useful in developing new  model-based,  probe-specific methods of    extracting target concentrations from  array fluorescence readings. \end{abstract}


\section{Introduction}

High density oligonucleotide microarrays are a widely used modern bio-technology tool  enabling the simultaneous 
testing for the presence as well as quantification of large numbers of genes in prepared target RNA samples. For a general  introduction to this technology we refer to the  celebrated paper    \cite{sinclair}  or  to \cite{arrays} for a more recent overview.  Among several  competing types of oligonucleotide microarrays, the Affymetrix GeneChip design appears to  be currently  one of the most common. GeneChip arrays consist of a substrate onto which short single strand DNA oligonucleotide probes have been synthesized using a photolithographic process. A 
chip surface is divided into some hundreds of thousands of regions typically tens of microns in size, with the DNA probes within each region being synthesized to a 
specific nucleotide sequence. The target RNA 
sample is hybridized onto the chip to form probe-target duplexes, and the chip is scanned to 
obtain fluorescence intensity readings from dyes incorporated during the laboratory procedures.
In principle, with suitable calibration, intensity readings are intended as a 'proxy'  measure of the concentration of matching target RNA in the 
sample.  However, due to optical noise, nonspecific hybridization, probe-specific effects, and measurement error, the  empirical measures of expression (i.e., the scanner-measured fluorescences) that summarize probe intensities  can often lead to imprecise and inaccurate results (see, e.g., \citet{p177}). 

It seems that some potentially significant improvement in relating the scanner readings of the probe intensities to the target genes concentrations could be obtained by using  a model-based approach   accounting for the physical processes driving hybridization. 
 Recently,  some authors  have begun to address this issues  by appealing to the dynamic adsorption models  well known in  physical chemistry literature (see, \citet{held}  or \citet{burden}). Such models  stemming from the physics of the chemical reactions involved are especially valuable  as  they could also help us in 
better  understanding of the physical processes driving hybridization 
and  lead to improvements in both microarrays design and performance. 

One of the most popular adsorption models considered in the context of microarrays (cf. e.g.,  \citet{hekstra} or  \citet{burden}) is the so-called {\em Langmuir model} (see next section) which in its simplest deterministic form  describes  the relationship between concentration and fluorescence levels of probe-target complexes by means of a  hyperbolic function.  In the context of microarrays  (in particular, GeneChips)  in order to   properly  account for the effects of multiple simultaneous hybridizations as well as the cross-hybridization due to competition  between  similarly sequenced targets for the same probe regions,  its seems that the stochastic version of the Langmuir model is  needed. The analysis of such a model was  carried out recently for instance in \citet{burden} or earlier in \citet{newt2, newt04} 
by means of adopting  the general results of \citet{denipati} on the fluctuations of the stochastic diffusion  equations around their stable equilibrium points.

The model for the stochastic  fluctuations of the  equation described by Dennis and Patil was cast  as a  boundary-free  problem and intended to provide a  continuous diffusion-type  approximation to the behavior of large biological systems as typically encountered in population dynamics problems. With no natural boundary restrictions  it was argued  in \citet{denipati} that the fluctuations around stable equilibria are approximately distributed as a gamma random variable. 
Based on this argument the gamma model for gene expressions  was since adopted by several authors in the context of analyzing microarray data (cf. e.g., \citet{newt2,newt,burden}). 

The simple extension of the Dennis and Patil results to microarrays setting, albeit appealing, seems to require further justification since the  microarray hybridization  models are neither continuous nor  boundary-free.   Whereas the continuous approximation to  the large discrete system seems easily justifiable,  it is not entirely clear what discrete system is being approximated by the boundary-free  diffusion  model (see \eqref{slang} below).  

The purpose of the current paper is to formally derive some simple closed-form  stochastic laws  approximating the equilibrium distributions of the discrete stochastic hybridization reactions under the explicit assumptions on the random noise terms which are consistent  with the stochastic Langmuir model but, unlike the latter, are not boundary-free.  The idea for the derivation is very simple. We start by noticing that the reaction rate equation of the deterministic Langmuir model may be considered as the usual approximation to the set of two coupled stochastic chemical reactions  on the finite state space (i.e., probe region size for GeneChips).  We then add  the stochastic forcing of  the Langmuir equation   as an   additional term to the original     birth and death rates of this discrete chemical system.
It turns out that the analysis of the equilibrium distribution of  this adjusted  birth and death process (we refer  to it below as the Langmuir BD process) when the number of states is large leads to the stochastic laws which are, under most circumstances,     consistent both with \citet{denipati} approach as well as \citet{burden} results.  However, due to the fact that we had based our analysis on a  finite discrete system, unlike previous approaches, ours gives more insight into the boundary behavior of the underlying discrete stationary process and its approximations. 
In particular,  when considering the discrete system it become obvious  that an  adjustment for   the saturation effect is needed  in the form  of  a Dirac-delta probability distribution at the boundary of the state space. This leads to an interesting consequence that the limiting stochastic law is not absolutely continuous as in Dennis and Patil result but rather has an   atom at the state space boundary.   We give some formal details of these findings in Section~3 below. 
 
Beyond the current introductory section  the paper is organized as follows. In the next section (Section~2) we offer  a brief overview  of some of the results related to the classical Langmuir adsorption-desorption  model  in our context. Our main theorem on the limiting stochastic law for the stationary distributions of the Langmuir birth-and-death process in large state space  along with some discussion is  presented in Section~3 along with an outline of the proof.  We conclude with some final remarks in Section~4.

\section{The Langmuir Model}   

In 1916 Irving Langmuir devised a simple model 
involving a thermodynamic equilibrium to predict the fraction of solid surface covered by an 
adsorbate as a function of its gas pressure \citep{langm}. 
The model  was later extended to liquid systems, where the equilibrium involved concentrations 
in solution.  In the Langmuir  model adsorbate and solvent molecules compete to adsorb on sites
 on the  surface of the powder and each site must be occupied by either a solvent molecule or an
  adsorbate molecule. For
 the hybridization reaction in 
 oligonucleotide DNA-microarrays  the same principle is applied 
in order to represent  competing adsorption and desorption of RNA molecules to form probe-target complexes at the chip 
surface (see, e.g., \citet{forman}).  

Let $\t=\t(t)\in (0,1)$ be the fraction of sites within a probe 
region occupied by probe-target complexes at time $t$ after the commencement of 
hybridization, and $d_1$ and $d_2$ be the forward adsorption and backward desorption rate constants respectively. The forward adsorption reaction is assumed 
to occur at a rate $d_1x(1-\t)$, proportional to the RNA-target concentration $x$ and 
fraction $(1-\t)$ of unoccupied probe sites. The backward  reaction  (desorption) 
is assumed to occur at a rate $d_2\,\t$, proportional to the fraction of occupied 
probe sites. In deterministic setting, the fraction of probe sites occupied by probe-target complexes is 
then given by the  reaction rate equation  known as the {\em Langmuir equation}
\begin{equation}\label{lang}
\frac{d\t}{dt}=d_1x-(d_1x+d_2)\t.
\end{equation}
The corresponding  equation incorporating the stochastic noise associated with both target and non-target specific  hybridization 
is given by the following stochastic version of \eqref{lang} herein refered to as the {\em stochastic Langmuir equation}. It has the form 
\begin{equation}\label{slang}
\frac{d\t}{dt}=d_1x-(d_1x+d_2)\t+\sqrt{g(\t)}\,Z_t,
\end{equation}
where $g(\t)\ge0$ is a known function of $\t$ and $Z_t$ is a Gaussian white noise process with unit variance.
The model described by the above equation is known in the literature as the {\em stochastic Langmuir model} and is a special case of a  diffusion model   considered e.g., by \citet{denipati}  in their study of stochastic fluctuations of populations about their stable equilibria.  We note that  in the present context of the Langmuir adsorption-desorption model,  the solution of the stochastic equation \eqref{slang}   is no longer bounded and  thus  \eqref{slang}  suffers  an obvious drawback   in the fact that   the function  $u$   has no physical interpretation for  $u>1$.  

We also note that \eqref{slang} 
 is concerned with  a single DNA (oligo) probe only, with the  effect of other probes replaced by a random noise term (stochastic forcing). In the context of modeling  a reaction network of simultaneous hybridization reactions on a DNA-microarray  \eqref{slang} is  therefore one of the simplest examples   of a  complex system  {\em decoupling} and {\em stochastic excitation}  (see, e.g., \citet{majda03,katsou}).
 
 In general, under mild regularity conditions on $g$, one can argue that the stationary solution of \eqref{slang} is approximately distributed as a gamma random variable around the deterministic system {\em  steady state} \citep{denipati}. However,  as noted by \citet{burden} for a linear choice of $g$, namely
\begin{equation}\label{g} g(\t)=Cd_1x\t\end{equation}with $C>0$ the equation  \eqref{slang} has  an {\em exact}  stationary gamma solution. That fact may be  easily inferred from the corresponding   Fokker-Planck equation (see, e.g., the monographs by \citet {vankamp} or \citet{EK86} and  the  references therein) which written in terms of the density of  $\t$, say  $\vp(\t,t)$,   is given by\footnote{Herein  we are primarily concerned with modeling the {\em internal}  fluctuations  of the stochastic system modeled by \eqref{slang}. Accordingly,  we interpret the stochastic equation \eqref{slang} in the sense of It\^o calculus. For the discussion of an alternative Fokker Planck equation \eqref{fp} using the Stratonovich calculus, see e.g., \citet{denipati} or, for more detail, \citet{vankamp}.}  \begin{equation}\label{fp}
\frac{\pa \vp(\t,t)}{\pa t}=\frac{\pa}{\pa \t}\left\{[(d_1x+d_2)\t-d_1x]\vp(\t,t)+\frac{C\,d_1x}{2}\,\frac{\pa[ \t\, \vp(\t,t)]}{\pa \t} \right\}.
\end{equation}  Solving for the steady state density, say   $\vp_0(u)$, gives \begin{equation}\label{bur_gam}
\vp_0(\t)\propto \t^{2/C-1}\exp\left(-2\,\frac{d_1x+d_2}{C\,d_1x }\,\t\right) \quad \text{for } \t\in [0,\infty)
\end{equation}
and zero otherwise.  

From the above considerations we see that adopting the stochastic model \eqref{slang} with no additional assumptions  may result in a stationary solution $\vp_0(u)$ being an absolutely continuous  distribution with positive support on  $(0,\infty)$.  This finding is, however,  not consistent with the experimental data which suggests that, at least for some values of the parameters, the saturated state $\t=1$ should  have  positive probability.  We also note an apparent  lack of  physical interpretation for the values $\t>1$ in the context of the original Langmuir model. 

An alternative approach to  modeling the dynamics of  hybridization  (or absorption-desorption) reaction is to analyze directly an underlying  discrete stochastic system which   \eqref{slang} intends to approximate.  We note that in our setting we have a simple one dimensional BD process described by one chemical species $Cmpx$ i.e., the amount of probe-target complex or, in other words, the number of occupied nucleotides in the probe region.   Hence, we consider  a system of two coupled chemical reactions  
 {\begin{eqnarray}\label{eq:sys}
\emptyset &\stackrel{b(\cdot)}\longrightarrow &Cmpx\nonumber\\
Cmpx &\stackrel{d(\cdot)}\longrightarrow&\emptyset
\end{eqnarray}} where $b(\cdot)$ and $d(\cdot)$ are  system-state dependent birth and death rates, respectively.  
In order to describe our approach  we need to specify    the form of these rate functions. To this end we shall define 
  a discrete, finite-state  version of the stochastic   Langmuir  adsorption-desorption.   

\begin{defn}[LBD Process]\label{def:lbd}{\em Langmuir BD process} is any  BD process with the set of states  $\{0,\ldots,N \}$   and the  birth and death  rates of the form  \begin{eqnarray*}
b(k) =&c_1(N-k)+\C(k,N)  \\
d(k) =&c_2\,k+\C(k,N)\end{eqnarray*} for $k=0,\ldots,N$.
Here  $c_1, c_2 $ are  some positive  real constants and the function $\C(\cdot,N)$ is assumed   to be of the form \begin{equation}\label{eq:nf}
\C(k,N)=\frac{N^2}{2}\,g(k/N)
\end{equation}
for $0<k<N$ and to satisfy the boundary conditions  ensuring the finiteness of the system space, i.e.  $\C(0,N)=\C(N,N)=0$. \end{defn}

In the LBD process the terms $c_1(N-k)$ and $c_2\,k$ are linear rates of birth and death as suggested by the deterministic Langmuir model \eqref{lang}.  The additional term $\C(\cdot,N)$ introduced into $b(k)$ and $d(k)$ is 
intended to model the noise of the {\em non-target} adsorption and desorption. For instance  on a GineChip  $\C(\cdot,N)$  accounts for the competition for the same RNA targets between different probe regions  with similar nucleotide sequences.  Assumption \eqref{eq:nf}  implies the ``density-dependent" form for the rates $b(\cdot), d(\cdot)$ (see, e.g., \citep{EK86} chapter 11) with the noise term $\C(\cdot,N)$ being of  higher order than the terms  $c_1(N-k)$ and $c_2\,k$.  

Note that  the (infinite) BD process with  the boundary-free rates (i.e., $d(k)$ and $b(k)$ given  as in Definition~\ref{def:lbd} but without requiring that $\C(0,N)=\C(N,N)=0$) may be approximated by the  solution of the Langmuir stochastic equation \eqref{slang} for large $N$. This may be informally  argued  as follows. Let $k_t$ be the state of the system \eqref{eq:sys} at $t\ge 0$, described as a difference  of two independent unit Poisson processes, say  $Y_1, Y_{-1}$, with random time changes (see, e.g., \citep{EK86} chapter 6) 
\begin{equation}\label{eq:see}
k_t=k_0+Y_{1}\left(\int_0^t b(k_s) ds\right) -Y_{-1}\left(\int_0^t d(k_s)ds\right).
\end{equation}
Since for any unit Poisson process $Y$ and large $N$ we have $ N^{-1/2}(Y(Nv)-Nv)\approx W(v)$  for any real $v$ with  $W(v)$ being the standard Brownian motion (SBM), the Poisson processes $Y_1, Y_{-1}$ may be approximated for large $N$ by independent SBM processes, say  $W_1,W_{-1}$. Denoting $u(t)=k_t/N$ 
this diffusion approximation of \eqref{eq:see} is 
\begin{align*}\label{eq:see2}
u(t)=& N^{-1}k_0+N^{-1/2} W_1\left(\int_0^t[c_1(1-u(s))+\frac{N}{2} g(u(s)) ]ds \right)\\ &-N^{-1/2}W_{-1}\left(\int_0^t[c_2u(s)+\frac{N}{2} g(u(s)) ]ds \right)+\int_0^t [c_1(1-u(s))+c_2u(s)]ds.
\end{align*}  which is distributionally equivalent to 
\begin{equation}\label{int_slang}
u(t)=\int_0^t [c_1(1-u(s))+c_2u(s)]ds +\int_0^t\sqrt{g(u(s))}dW(s)+o_P(1)
\end{equation}
and hence in the limit to the integral version of \eqref{slang}  (see \citet{EK86} chapter 11 for details).
 
Of course,  depending on the form of  $\C(\cdot,N)$ we shall have different forms of the LBD process.  In order to cast  our results  somewhat parallel to the model \eqref{slang} under  a linear form of $g$ in \eqref{g},   herein we consider only  $\C(k,N)$ given by the functions $\C_1,\C_2,\C_3$ defined below, with the  corresponding models  henceforth  referred to  as  $M_1, M_2$,  $M_3$, respectively. 
\begin{align*}
 \C_1(k,N) &=c_3Nk \qquad\qquad\text{ for } 0\le k<N \text{ and }\C_1(N,N)=0 &(M_1)\\
\C_2(k,N) &=c_3N(N-k)\quad\text{ for }  0<k\le N \text{ and } \C_2(0,N)=0 &(M_2)\\
\C_3(k,N)&=c_3k(N-k)\ \quad\text{ for }  0\le k\le N &(M_3)
\end{align*}

Note that if we disregard the boundary condition  $\C_1(N,N)=0$ then the  model $M_1$ is a discrete  analogue of \eqref{slang} with  the  choice of $g$ given by \eqref{g} in the sense that BD process \eqref{eq:sys} (or \eqref{eq:see}) may be approximated by \eqref{int_slang} for large $N$,  leading  to the  Fokker-Planck equation given in \eqref{fp}  and, consequently, to \eqref{bur_gam} with  
\begin{eqnarray}\label{cd}
c_1 & = & d_1x \nonumber\\
c_2 & = & d_2 \\
c_3  & = & Cd_1x/2.\nonumber 
\end{eqnarray}
This casting of the equation \eqref{slang} as an approximation to $M_1$ gives also  additional  insight into the somewhat mysterious choice of the form of function $g$ in \eqref{g}.
Considering  the rates in  $M_1$ it becomes clear that the choice of \eqref{g} is 
a reflection of two  implicit assumptions concerning the microarray  hybridization reactions, namely that  (i)  the  level of the target-specific signal in the probe region has  lower magnitude than  the level of non-specific signal (i.e., signal noise)  and (ii) the non-specific signal noise is proportional to  the total system (i.e., probe region) size  as well and the current system state and the target concentration. 

Note that the model  $M_2$ is simply a  'reflection' of  $M_1$ obtained by considering the amount of unoccupied probe region $N-Cmpx$ instead of the amount of $Cmpx$. Model $M_2$ is thus not concerned with saturation but rather with {\em threshold effect} of the probe adsorption.  This phenomena occurs when the LBD process attains an {\em empty} state with positive probability.  

Note also that  both $M_1$ and $M_2$ rate functions have discontinuities at the boundary. This is in contrast with   the model $M_3$ which enjoys   smooth  boundary conditions  with no discontinuities. 
In general such discontinuities in rate functions for BD processes 
prevent the direct application of an approximation of the form \eqref{int_slang}, however it turns out that for  an LBD processes $M_1-M_3$ their   stationary distributions may be approximated more directly. 

\section{Limit Theorem}
In this section we state  and prove the  main result of the paper, namely the limit theorem for the stationary distributions of the LBD processes under the models $M_1-M_3$. The proof we give herein is quite elementary  and  is  based on the fact that for  one dimensional birth-and-death processes with  bounded state space and polynomial rates, the moments of their  limiting distributions must be uniquely determined by the corresponding detailed balance (reversibility) conditions.  At this point it is perhaps also worth noticing that even though herein we have restrict ourselves only to   the models $M_1-M_3$  as the  most relevant for  the type of asymptotic behavior described by \eqref{slang} and \eqref{g}, it is not difficult to see that the method of the proof allows one to easily extend the result to   any LBD processes   with polynomial-type birth and death rates. 
This,  at least in principle, allows then to obtain limit theorems for  the discrete versions of \eqref{slang} with any  function $g$ continuous on (0,1) and  continuously extendable to [0,1] where it may  be always uniformly approximated by    polynomials. However,  such considerations go beyond our current scope.
 
In order to state the theorem we shall need some additional notation.  For $z,\gamma>0$ denote the incomplete gamma function 
by $\G(z,\g)=\int_0^{\g} s^{z-1} \exp(-s)\,ds $ and for any $\a,\b>0$  let $IG(\a,\b,1)$ denote an incomplete gamma random variable with the density function  $f_{\a,\b}(x)=\G(\a,\b)^{-1}\, \b^\a\,x^{\a-1}\,\exp(-x\b)$ for $x\in (0,1)$ and zero otherwise.  Let us denote by $F_{\a,\b}$ the distribution function of  $IG(\a,\b,1)$. We introduce  the following definition.

\begin{defn}[LIG Distribution]\label{def:lig} 
We say that the random variable  has the {\em Langmuir-incomplete gamma} (LIG)  distribution with parameters $\a,\b$ satisfying $\b>\a>0$ if  its distribution function  is given by the mixture 
$$G=(1-\pi_{\a,\b}) F_{\a,\b}+\pi_{\a,\b}\,\delta_1$$ where  $F_{\a,\b}$ is the distribution function for  $IG(\a,\b,1)$ random variable, $\delta_1$ is the distribution function of a degenerate  random variable with mass concentrated at one and   \begin{equation}\label{eq:pidef}
\pi_{\a,\b}=\frac{\b^\a}{\b^\a+\G(\a,\b)\,(\b-\a)\exp(\b)}.\end{equation} 
\end{defn}
Below we denote the Langmuir incomplete gamma distribution with parameters $\a,\b$ by $LIG(\a,\b)$. We shall also denote by  $Beta(\a,\b)$ the usual beta distribution with  parameters $\a,\b>0$ and the density  $h(x)=\G(\a+\b)\G(\a)^{-1}\G(\b)^{-1}\int_0^1x^{\a-1}\,(1-x)^{\b-1}\,dx$ for $x\in(0,1)$ and zero otherwise. We have the following

\begin{thm}[Limit Theorem for a Stationary Distribution of an LBD Process]\label{thm:1} Let $X_N^{(i)}$ be the stationary distributions  of LBD Process  $M_i$ for $i=1,2,3,$ and let    $a=c_1/c_3$ and $b=(c_1+c_2)/c_3$,  as well as $Y^{(i)}_N=  X_N^{(i)}/N$. Then, as $N\to \infty$ we have the weak convergence 
$$Y^{(i)}_N\stackrel{D}\rightarrow {\cal Z}_i\qquad i=1,2,3$$ where the limiting random variables ${\cal Z}_i$ are as follows
\begin{itemize} 
\item[(i)] $ {\cal Z}_1$ is  $LIG(a,b)$, 
\item[(ii)] ${\cal Z}_2 $ is such that $1-{\cal Z}_2 $ is $LIG(b-a,b)$,
\item[(iii)] ${\cal Z}_3$ is  $Beta (a,b-a).$
\end{itemize} 
\end{thm}

Before  discussing the proof of this result, some  remarks  are perhaps in order. 

\begin{itemize}
\item As  it shall  become clear from the proof, it turns out that  all any    LBD processe (hence, also $M_1-M_3$ )  have the correct Langmuir mean  given by the stationary solution of the deterministic equation \eqref{lang} with  the  $c_1$ and $c_2$ constants as  in   \eqref{cd}.   In that sense an LBD process  may be viewed as a discrete analogy of  continuous models \eqref{lang} and \eqref{slang}  with specific   functions $g$ in the latter   related to  LBD   via \eqref{eq:nf}.  One should stress, however, the fundamental difference between the approach  to approximating a stochastic 
equilibrium of a discrete system \eqref{eq:sys} offered by Theorem~\ref{thm:1}  and that based on  analyzing the  equilibrium distribution of the diffusion approximation   \eqref{slang} outlined in  \eqref{eq:see} and \eqref{int_slang}.  In   Theorem~\ref{thm:1} one considers a sequence of  stationary distributions of  LBD processes indexed by the size of the state space $N$ and derives its limit as $N$ increases. The approximation  via \eqref{slang} is based on approximating the entire discrete process (not just its equilibrium distribution) for large $N$ and then  deriving a stationary distribution of the approximation. 

\item Despite the very different models behind them,  if $\pi_{a,b}\approx 0$ then   $(i)$ of Theoem~1 specializes to the  result  on the stationary density \eqref{bur_gam} obtained via  \eqref{fp}. Thus when $\pi_{a,b}\approx 0$ our theorem fomally justifies the use of gamma approximation for modeling hybridization reactions under the boundary-free model described by \eqref{slang} and the assumption \eqref{g}.  \item If the condition $\pi_{a,b}\approx 0$ is not satisfied then there may be a significant difference between the stationary distribution obtained from the boundary-free type analysis via the stochastic equation \eqref{slang} and the LBD process analysis. This is due to the fact that the LBD analysis takes properly  into account the discontinuities in the   rate functions whereas the continuous model \eqref{slang} does not.   

\item The theorem indicates that both $M_1$ and $M_2$ models which incorporate the linear noise term $\C(\cdot,N)$ into their rate functions are amenable  to the gamma-type approximation of  their stationary distributions perhaps after some adjustment for  the  bounded  state-space. In contrast,  the  model $M_3$  with the quadratic and 'boundary symmetric' noise  term   yields a different type of stationary distribution (i.e, beta) with no boundary effects.
\end{itemize}

 It seems that there are several ways of arriving at the result of the theorem. Herein  we have chosen the method of the proof which is perhaps slightly convoluted but on the other hand almost completely elementary and thus fully accessible to readers without  extensive background in stochastic processes theory. 


In order to provide a proof of Theorem~1 we shall need two auxiliary  results stated below as Lemmas~1 and 2. The first one of them concerns some elementary properties of  the moments of a  LIG distribution.

\begin{lem}\label{lem:1}  Let $Z$ be a random variable distributed according to $LIG(\a,\b)$. Then for any integer $m\ge 0$ we have 
\begin{equation}
\label{eq:lig_recur}
EZ^{m+1}=\frac{m+\a}{\b}EZ^{m}-\frac{m}{\b}\pi_{\a,\b}
\end{equation}

\end{lem}

\begin{proof}[Proof of Lemma~\ref{lem:1}]
 Let $W$ be a random variable distributed according to $IG(\a,\b,1)$. Elementary calculation based on the integration by parts shows that for any integer $m\ge 0$ 
\begin{equation*}
\label{eq:1 }
EW^{m+1}=\frac{m+\a}{\b}EW^{m}-\frac{\,\b^{\a-1}\,\exp(-\b)}{\G(\a,\b)}.
\end{equation*}
In view of the above and by the definition of $Z$ we have for any integer $m\ge 0$
\begin{align*}  E Z^{m+1} & = (1-\pi_{\a,\b})\,EW^{m+1}+\pi_{\a,\b} \\
& = (1-\pi_{\a,\b}) \left[\frac{m+\a}{\b}EW^{m}-\frac{\,\b^{\a-1}\,\exp(-\b)}{\G(\a,\b)} \right]+\pi_{\a,\b} \\
& = \frac{m+\a}{\b}  \left[(1-\pi_{\a,\b})\,EW^m +\pi_{\a,\b} \right]+\left[1-\frac{m+\a}{\b} \right]\,\pi_{\a,\b} -(1-\pi_{\a,\b})\,\frac{\b^{\a-1}\,\exp(-\b)}{\G(\a,\b)}\\
& = \frac{m+\a}{\b}  \left[(1-\pi_{\a,\b})\,EW^m +\pi_{\a,\b} \right]-\frac{m}{\b}\pi_{\a,\b} -\frac{\b-\a}{\b} \pi_{\a,\b} -(1-\pi_{\a,\b})\,\frac{\b^{\a-1}\,\exp(-\b)}{\G(\a,\b)} \\
& = \frac{m+\a}{\b}  \,EZ^m -\frac{m}{\b}\pi_{\a,\b}.  \end{align*}
\end{proof}
Our second lemma is as follows.
\begin{lem}\label{lem:2} Let $\a,\b>0$ be arbitrary. Consider $N\to \infty$. For  any non-increasing real sequence $\a_N \downarrow \a>0$ satisfying $ (\a_N-\a)\,\log N \to 0$  and any  real sequence $\b_N\to  \b>0$   we have 
$$ N^{-\a_N}\,\sum_{k=0}^N \frac{\G(\a_N+k)}{k!} \left(1-\frac{\b_N}{N}\right)^{k}\to\int_{0}^1 x^{\a-1}\,e^{-\b x}\,dx=\b^{-\a}\,\G(\a,\b).$$
\end{lem}
\begin{proof}[Proof of Lemma~\ref{lem:2}] Assume first that $\a_N\equiv \a$.
For given $\a,\b>0$ define $k(N)=[\d\log N]$ where $\d$ is a fixed positive number such that $\d<\a$ and  $[x]$ denotes  the largest integer not greater than $x$.  Write \begin{align*}
N^{-\a}\,\sum_{k=0}^N \frac{\G(\a+k)}{k!} \left(1-\frac{\b_N}{N}\right)^{k}&=N^{-\a}\,\sum_{k=0}^{k(N)} \frac{\G(\a+k)}{k!} \left(1-\frac{\b_N}{N}\right)^{k}+N^{-\a}\,\sum_{k=k(N)+1}^N \frac{\G(\a+k)}{k!} \left(1-\frac{\b_N}{N}\right)^{k}\non \\
&=(I)+(II)\end{align*}
We first show $(I)\to 0$ as $N\to \infty$. To this end note
\begin{align*}
(I)&\le N^{-\a}\,\sum_{k=0}^{k(N)} \frac{\G(\a+k)}{k!}\le N^{-\a}\,\sum_{k=0}^{k(N)} \frac{(\a+k)^k}{k!}\\
& \le N^{-\a}\,\sum_{k=0}^{k(N)} \frac{(\a+k(N))^k}{k!}\le N^{-\a}e^{(\a+k(N))}\to 0 \quad \text{as}\ N\to \infty. 
\end{align*}
Now we argue that \begin{equation}\label{eq:1}(II)\to \b^{-\a}\G(\a,\b) \quad\text{as}\quad N\to \infty. \end{equation} To this end, recall  the following version of the Gauss formula \begin{equation}\label{eq:GF}
\frac{\G(\a+k)}{k!\,k^{\a-1}}\to 1 \quad\text{as}\quad k\to \infty. \end{equation} In view of the above it follows that for any given $\ve\in(0,1)$ and $N$ sufficiently  large we  have 
$$  (1-\ve)  N^{-\a}\hspace{-10pt}\sum_{k=k(N)+1}^{N} k^{\a-1}e^{-k\b/N} \le (II)\le (1+\ve)  N^{-\a}\hspace{-10pt}\sum_{k=k(N)+1}^{N} k^{\a-1}e^{-k\b/N}.$$ Since the expression  $ N^{-\a}\sum_{k=k(N)+1}^{N} k^{\a-1}e^{-k\b/N}= N^{-1}\sum_{k=k(N)+1}^{N} (k/N)^{\a-1}e^{-k\b/N}$ is seen to be the Riemann sum for $\b^{-\a}\G(\a,\b)$ taking   $N\to \infty$ gives 
$$(1-\ve)\b^{-\a}\G(\a,\b)\le \lim_N\, (II)\le(1+\ve)\b^{-\a}\G(\a,\b). $$ Since $\ve$ may be taken arbitrarily close to zero,
the  relation \eqref{eq:1} follows and yields the  assertion  of the lemma with $\a_N\equiv \a$. To complete the proof for an arbitrary  sequence $\a_N$  note that we need in essence  only to argue that  $N^{\a_N-\a}\to 1$  as $N\to \infty$ (this follows by assumption) and that the formula \eqref{eq:GF} holds with $\a$ replaced by $\a_k$ (since $\a_k$ is monotone).  By the continuity of gamma function,
$\G(\a)/\G(\a_N)\to 1$ as $N\to \infty$.  This  and   \eqref{eq:GF} entail $$\frac{\G(\a_k)\,\a\,(\a+1)\cdots(\a+k-1)}{k!\,k^{\a-1}}\to 1 \qquad \text{as}\  k\to \infty.$$ The relationship \eqref{eq:GF} with $\a$ replaced by $\a_k$  will now follow  if we can argue 
that \begin{equation} \label{eq:prod2one}\prod_{s=0}^k \frac{\a+s}{\a_k+s}\to 1\end{equation} as $k\to \infty$. To this end note that $$\log \left(\prod_{s=0}^k \frac{\a+s}{\a_k+s}\right) =\sum_{s=0}^k \log\left( \frac{\a_k+s}{\a+s}\right)\le (\a_k-\a) \sum_{s=0}^k \frac{1}{s+\a}\le 2\,(\a_k-\a) \log k\to 0 $$ by our assumption on $\a_k$ and thus \eqref{eq:prod2one} follows. This, however, yields the assertion  of the lemma, since for  sufficiently large $N$ 
$$ N^{-\a}\! \sum_{k=k(N)}^N \frac{\G(\a+k)}{k!} \left(1-\frac{\b_N}{N}\right)^{k} \le N^{-\a}\! \sum_{k=k(N)}^N \frac{\G(\a_N+k)}{k!} \left(1-\frac{\b_N}{N}\right)^{k} \le N^{-\a}\! \sum_{k=k(N)}^N \frac{\G(\a_k+k)}{k!} \left(1-\frac{\b_N}{N}\right)^{k} $$ and we have just shown that the first and the last of the expressions above tend to $\b^{-\a}\,\G(\a,\b)$ as $N\to \infty$. 
  \end{proof}
Having established the assertions of the lemmas above, we are finally in a position to   prove the  result given in Theorem~1. 
\begin{proof}[{\bf Proof of Theorem~\ref{thm:1}}. Part (i)]
 Denote  by $X$ the random variable $X_N^{(1)}$ and set $P(X=k)=p(k)$ for $k=0\ldots,N$. Let $m\ge 0$ be an integer.   Multiplying by $k^m$ the detailed balance equation 
\begin{equation}\label{eq:bal}
p(k+1)\,d(k+1)=p(k)\,b(k)
\end{equation}
and then summing over $k=0,\ldots, N-1$ we obtain under $M_1$ model 
\begin{equation*}
(c_2+c_3\,N)\,\sum_{k=0}^{N-1}\,k^m\,(k+1)\,p(k+1)-c_3\,(N-1)^l\,p(N)\,N^2=\sum_{k=0}^{N-1}\,k^m\left[p(N)(c_1(N-k)+c_3\,N\,k) \right].
\end{equation*}
Expanding now $k^m=(k+1-1)^m$ on the right hand side we get \begin{align*}
&(c_2+c_3\,N)\, \sum_{k=0}^{N-1}\sum_{s=0}^m \binom{m}{s}(-1)^{m-s}\,(k+1)^{s+1}\,p(k+1)\ - c_3\,(N-1)^m\,N^2\,p(N)\\
&= \sum_{k=0}^N [c_1(N-k)+c_3Nk]\,k^m\,p(k)-c_3\,N^{m+2}p(N).
\end{align*}
Denoting $\wm_m(N)=EX^m$ we may rewrite the above relationship as 
\begin{align*}
&(c_2+c_3\,N)\,\sum_{s=0}^{m}\binom{m}{s}(-1)^{m-s}\wm_{m+1}(N)-c_3\,(N-1)^m\,N^2\,p(N)\\
&= c_1\,N\wm_m(N)+(c_3\,N-c_1)\wm_{m+1}(N)-c_3\,N^{m+2}\,p(N).
\end{align*}
We  set  $\mu_m=\lim_{N\to \infty}\,\wm_m(N)/N^m $. Note that   dividing both sides by $N^{m+1}$  and taking $N\to \infty$ gives 
\begin{equation*}
(c_1+c_2)\mu_{m+1}=(m\,c_3+c_1)\,\mu_m -c_3\,m\,\p
\end{equation*}
where \begin{equation}
\label{eq:pstar}
\p=\lim_{N\to \infty}p(N).
\end{equation}

\noi Assuming for a moment  that $\p$ exists  and is finite (this follows from \eqref{eq:pstar2} below) we see that  these considerations give the following recursive relationship for the limiting moments of $Y^{(1)}_N$ 
$$\mu_{m+1}=\frac{m+c_1/c_3}{(c_1+c_2)/c_3}\,\mu_m-\frac{m\,\p}{(c_1+c_2)/c_3}\qquad\text{for}\ m=0,1\ldots.
$$
or in terms of $a,b$
\begin{equation}\label{eq:recur1} \mu_{m+1}=\frac{m+a}{b}\,\mu_m-\frac{m}{b}\p\qquad\text{for}\ m=0,1\ldots.\end{equation}
We note that in view of  $\mu_0=1$  the solution to the above recursive equation is unique. Moreover, we note that since the  support of the sequence  $\{Y^{(1)}_N\}_{N=1}^\infty$
is contained in  the closed interval [0,1]  then (i) the corresponding sequence of probability measures   is tight and (ii) any of its weak limits must be  a probability measure whose  moments satisfy \eqref{eq:recur1}.  Since the probability distributions on bounded intervals are  uniquely determined by their  moments  it follows that  as $N\to \infty$
\begin{equation}\label{eq:conv1}
Y^{(1)}_N\stackrel{D}\rightarrow \Z_1 
\end{equation}
for some random variable $\Z_1$ with moments $\mu_m$ given by \eqref{eq:recur1}. To complete the proof of part $(i)$ we need to show only that $\Z_1$ is a $LIG(a,b)$ random variable  as given in Definition~\ref{def:lig}. 
Since $\Z_1$ is identified completely by its moments, it suffices to show that the moments of the random variable $LIG(a,b)$ satisfy the recursive relation \eqref{eq:recur1}. This follows by  Lemma~1 provided that  
\begin{equation}\label{eq:pstar2}
p^\ast=\pi_{a,b}
\end{equation}  
where $\pi_{a,b}$ is given by \eqref{eq:pidef}.

In order to argue \eqref{eq:pstar2} we again consider the detailed balance equation \eqref{eq:bal} 
\begin{align}\label{eq:delta} p(N) & = p(N-1)\,b(N-1)/d(N)=\frac{p(N-1)\,b(N-1)/d(N)}{p(N-1)\,b(N-1)/d(N)+\sum_{k=0}^{N-1} p(k)}\nonumber\\
& = \frac{p(N-1)\,b(N-1)/d(N)}{p(N-1)\,b(N-1)/d(N)+\sum_{k=0}^{N-1} p(k)}=\frac{\Delta_N}{\Delta_N+1} \end{align} where we define $$\Delta_N= \frac{p(N-1)\,b(N-1)/d(N)}{\sum_{k=0}^{N-1} p(k)}.$$

We note that again by \eqref{eq:bal} we get under the model $M_1$ the following form of $\Delta_N$
\begin{align*}
\Delta_N & = \frac{c_2+c_3\,N}{c_2\,N!}\ \left( \frac{c_3\,N-c_1}{c_2+c_3\,N}\right)^N\  \frac{\prod_{s=0}^{N-1}(s+\frac{c_1\,N}{c_3\,N-c_1})}{\sum_{k=0}^{N-1} \frac{1}{k!}\left(\frac{c_3\,N-c_1}{c_2+c_3\,N}\right)^k\prod_{s=0}^{k-1}(s+\frac{c_1\,N}{c_3\,N-c_1})} \\
& =  \frac{c_2+c_3\,N}{c_2\,N!}\ \left(\frac{c_3\,N-c_1}{c_2+c_3\,N}\right)^N\  \frac{\G(N+\frac{c_1\,N}{c_3\,N-c_1})}{\sum_{k=0}^{N-1} \frac{1}{k!}\left(\frac{c_3\,N-c_1}{c_2+c_3\,N}\right)^k\G(k+\frac{c_1\,N}{c_3\,N-c_1})}.
\end{align*}
Denote $$a_N=\frac{c_1\,N}{c_3\,N-c_1}\qquad b_N=\frac{(c_1+c_2)\,N}{c_3\,N+c_2}\, $$ 
then 

\begin{align*}
\Delta_N & =\frac{(b-a+N)}{(b-a)\,N}\ \frac{\G(N+a_N)}{N!\, N^{a_N-1}}\, \frac{\left( 1-b_N/N\right)^N}{N^{-a_N}\,\sum_{k=0}^{N-1} \frac{1}{k!} \left(1-b_N/N\right)^{k} \G(k+a_N)}.
\end{align*} Applying now the Gauss formula \eqref{eq:GF} and using the result  of Lemma~2 we  conclude  that 
$$ \lim_N \Delta_N=\frac{b^a}{\exp(b)\,\G(a,b)}$$ and 
hence \eqref{eq:pstar2} follows  by \eqref{eq:delta} which completes the proof of part $(i)$ of the theorem. 
\smallskip

\noindent {\em Part (ii).} The result follows by applying part $(i)$ to the random variable $N-X_N^{(2)}$.

\smallskip

\noindent {\em Part (iii).} For the proof of the last part of the  theorem we denote  now by $X$   the random variable $X_N^{(3)}$ and otherwise retain the notation from part $(i)$. Multiplying  \eqref{eq:bal} by $k^m$ and summing as before, we obtain under $M_3$ $$\sum_{k=0}^N k^m\, p(k+1)[c_2(k+1)+c_3(k+1)(N-k-1)]=\sum_{k=0}^N k^m\,p(k)[c_1(N-k)+c_3k(N-k)].$$ 
The above, by an argument similar to the one used in $(i)$,  gives  the relationship
$$ (c_2+c_3N) \sum_{s=0}^{m} (-1)^{m-s}\binom{m}{s}\wm_{s+1}-c_3 \sum_{s=0}^{m} (-1)^{m-s}\binom{m}{s}\,\wm_{s+2} = c_1N\,\wm_m-c_1\wm_{m+1}+c_3N\,\wm_{m+1}-c_3\,\wm_{m+2}. $$
Denoting  $\mu_m=\lim_{N\to \infty}\,\wm_m(N)/N^m $, dividing both sides by $N^{m+1}$  and taking $N\to \infty$ gives  somewhat simpler then  \eqref{eq:recur1} recursion formula  for  the limiting moments, namely  
$$\mu_{m+1}=\frac{m+c_1/c_3}{m+(c_1+c_2)/c_3}\,\mu_m\qquad\text{for}\ m=0,1\ldots.
$$ Writing the above in    terms of $a,b$
$$\mu_{m+1}=\frac{m+a}{m+b}\,\mu_m\qquad\text{for}\ m=0,1\ldots
$$ we obtain the familiar relationship between the moments of $Beta(a,b-a)$ distribution. This, along with the  tightness of measures argument similar to the one used in $(i)$ completes the proof of part $(iii)$.  
\end{proof}

\section{Conclusions} 
Herein we have derived a  limit theorem for stationary distributions of some special  birth-and-death processes related to the Langmuir dynamic adsorption-desorption model. Such a model is of interest in the context of the microarray hybridization reactions if one may assume that the fluorescence signal on the array is approximately a  realization of a chemical  Langmuir equilibrium of  the adsorption and desorption reactions between the target  mRNA molecules and the  DNA probes.  Whereas this assumption may be questionable for  long (hundred basis or more)  probes,  it  seems reasonable for  the short ones, like e.g., the 25-mers used on  many  Affymetrix  chips.  Indeed, in the context of Affymetrix GeneChip arrays, the gamma-type  approximation to   the gene expression data based on an ad-hoc  Langmuir-like equilibria argument have been proposed in the literature as a way of enhancing the data analysis.  Our current result  gives a rigorous justification of  the use of  truncated-gamma and beta-type distributions in order  to approximate the fluorescence readings of the probe-RNA complexes obtained in course of an Affymetrix microarray experiment. It also    explains some experimentally observed behavior of these readings 
like  e.g. the signal saturation and the signal thresholding phenomena.  

The potential usefulness of our approximation results  stems also from  the fact that they allow one to describe the  theoretical means of measured fluorescence intensity readings  by  three parameter hyperbolic response functions which can be  obtained as  solutions of the corresponding deterministic Langmuir equations. In general, these response functions  for  specific probes  shall be only   sequence dependent  and could be  therefore used universally in  all experiments involving a particular probe sequence. 
Our results imply also  that the  fold changes in RNA-target concentration are  not linearly related to fold changes in fluorescence intensity readings, as is often generally assumed. 

  As pointed out by some authors (cf. \citet{p177})  the formidable challenge in microarray experiments is to  establish a reliable  algorithm for extracting the true RNA concentrations measurements from the probes fluorescence intensity readings.  We believe that the results of this paper could perhaps get us a step closer  to that goal. 
  
 \section*{Acknowledgement} 
 This research is partially sponsored by the National Science Foundation  under grant DMS0553701 as well as by 
the  Center for Environmental Genomics and Integrative Biology  at the University of Louisville which receives funding from the   National Institute of Environmental Health  Sciences under grant 1P30ES014443. 

  The authors wish to acknowledge an anonymous referee whose comments helped us improve the original manuscript. 
\bibliographystyle{apalike}
\bibliography{langmuir}

 \end{document}